\documentclass[12pt,reqno]{amsart}
\usepackage{amssymb}
\usepackage[latin1]{inputenc}
\usepackage{dsfont}
\usepackage{hyperref}
\usepackage{fullpage}
\newtheorem{theorem}{Theorem}[section]

\newtheorem{proposition}[theorem]{Proposition}

\theoremstyle{definition}
\newtheorem{definition}[theorem]{Definition}

\theoremstyle{remark}

\numberwithin{equation}{section}
\begin{document}
\title{A new construction of the $\sigma$-finite measures associated with submartingales of class $(\Sigma)$
(Une nouvelle construction des mesures $\sigma$-finies associ\'ees aux sous-martingales de classe
$(\Sigma)$)}
\author{Joseph Najnudel, Institut f\"ur Mathematik, Universit\"at Z\"urich, Winterthurerstrasse 190,
8057-Z\"urich, Switzerland, joseph.najnudel@math.uzh.ch}
\author{Ashkan Nikeghbali,
          Institut f\"ur Mathematik, Universit\"at Z\"urich, Winterthurerstrasse 190,
8057-Z\"urich, Switzerland, ashkan.nikeghbali@math.uzh.ch}
\date{\today}
\begin{abstract}
In \cite{NN1}, we prove that for any submartingale $(X_t)_{t \geq 0}$ of class $(\Sigma)$, defined
 on a filtered probability
space $(\Omega, \mathcal{F}, \mathbb{P}, (\mathcal{F}_t)_{t \geq 0})$, which satisfies some technical conditions, 
one can construct a $\sigma$-finite measure $\mathcal{Q}$ on $(\Omega, \mathcal{F})$, such that for all 
$t \geq 0$, and for all events $\Lambda_t
 \in \mathcal{F}_t$: $$ \mathcal{Q} [\Lambda_t, g\leq t] = \mathbb{E}_{\mathbb{P}} [\mathds{1}_{\Lambda_t} X_t]$$
where $g$ is the last hitting time of zero of the process $X$. Some particular cases of this construction 
  are related with Brownian penalisation or mathematical finance. In this note, we give a  
simpler construction of $\mathcal{Q}$, and we show that an analog of this measure can 
also be defined for discrete-time submartingales. \\
Dans \cite{NN1}, nous prouvons que pour toute sous-martingale $(X_t)_{t \geq 0}$ de classe $(\Sigma)$, d\'efinie
sur un espace de probabilit\'e filtr\'e $(\Omega, \mathcal{F}, \mathbb{P}, (\mathcal{F}_t)_{t \geq 0})$,
 satisfaisant certaines conditions techniques, 
on peut construire une mesure $\sigma$-finie $\mathcal{Q}$ sur $(\Omega, \mathcal{F})$, telle que pour tout
$t \geq 0$, et pour tout \'ev\'enement $\Lambda_t
 \in \mathcal{F}_t$: $$ \mathcal{Q} [\Lambda_t, g\leq t] = \mathbb{E}_{\mathbb{P}} [\mathds{1}_{\Lambda_t} X_t]$$
o\`u $g$ est le dernier z\'ero de $X$. Certains cas particuliers de cette construction sont li\'es aux 
p\'enalisations browniennes ou aux math\'ematiques financi\`eres. Dans cette note, nous donnons une 
construction plus simple de $\mathcal{Q}$, et nous montrons qu'un analogue de 
cette mesure peut aussi \^etre d\'efini pour des sous-martingales \`a temps discret.
\end{abstract}
\maketitle
\section{Version fran\c caise abr\'eg\'ee}
Les sous-martingales de classe $(\Sigma)$ ont \'et\'e introduites par Yor (voir \cite{Y}) et 
certaines de ses propri\'et\'es les plus importantes ont \'et\'e \'etudi\'ees par Nikeghbali dans \cite{N}. 
Ces sous-martingales sont en particulier li\'ees \`a certains probl\`emes de math\'ematiques financi\`eres
(voir \cite{MRY}, \cite{BeYor}, \cite{CNP}) et aux p\'enalisations browniennes (voir \cite{NRY}). 
La class $(\Sigma)$ est d\'efinie de la mani\`ere suivante:
 \begin{definition}
Soit $(\Omega,\mathcal{F},(\mathcal{F}_t)_{t \geq 0},\mathbb{P})$ un espace de probabilit\'e filtr\'e. Une 
sous-martingale $(X_t)_{t \geq 0}$ est de class $(\Sigma)$, si pour tout $t \geq 0$, $X_t = N_t + A_t$,
o\`u $(N_t)_{t \geq 0}$ et $(A_t)_{t \geq 0}$ sont des processus $(\mathcal{F}_t)_{t \geq 0}$-adapt\'es 
satisfaisant les conditions suivantes:
\begin{itemize}
\item $(N_t)_{t \geq 0}$ est une martingale c\`adl\`ag. 
\item $(A_t)_{t \geq 0}$ est un processus croissant, continu, tel que $A_0 = 0$;
\item La mesure $(dA_t)$ est port\'ee par l'ensemble $\{t \geq 0, X_t = 0 \}$.
\end{itemize}
\end{definition}
\noindent
Dans \cite{NN1}, nous prouvons qu'on peut associer une mesure $\sigma$-finie \`a toute sous-martingale 
de classe $(\Sigma)$, d\'efinie sur un espace de probabilit\'e filtr\'e 
$(\Omega,\mathcal{F},(\mathcal{F}_t)_{t \geq 0},\mathbb{P})$ satisfaisant une certaine condition 
technique, appel\'ee propri\'et\'e (NP). L'\'enonc\'e est le suivant:
\begin{theorem} \label{tout}
Soit $(X_t)_{t \geq 0}$ une sous-martingale de classe $(\Sigma)$ (en particulier $X_t$ est int\'egrable
pour tout $t \geq 0$), d\'efinie sur un espace de probabilit\'e filtr\'e 
$(\Omega, \mathcal{F}, \mathbb{P}, (\mathcal{F}_t)_{t \geq 0})$ qui satisfait la propri\'et\'e (NP).
En particulier, $(\mathcal{F}_t)_{t \geq 0}$
 est continue \`a droite et $\mathcal{F}$ est
la tribu engendr\'ee par $\mathcal{F}_t$, $t \geq 0$. 
Dans ces conditions, il existe une unique mesure $\sigma$-finie $\mathcal{Q}$, d\'efinie sur 
$(\Omega, \mathcal{F}, \mathbb{P})$ et telle que si $g = \sup\{t \geq 0, X_t = 0 \}$:
\begin{itemize}
\item $\mathcal{Q} [g = \infty] = 0$;
\item Pour tout $t \geq 0$, et pour toute variable al\'eatoire $F_t$ born\'ee, $\mathcal{F}_t$-mesurable,
$$\mathcal{Q} \left[ F_t \, \mathds{1}_{g \leq t} \right] = \mathbb{E}_{\mathbb{P}} \left[F_t X_t \right].$$
\end{itemize}
\noindent
\end{theorem} 
\noindent
La propri\'et\'e (NP) est d\'efinie pr\'ecis\'ement dans \cite{NN2}, et intervient en particulier 
en raison de probl\`emes recontr\'es lors d'extensions de familles compatibles de mesures de probabilit\'es
(ces probl\`emes sont \'egalement discut\'es par Bichteler dans \cite{Bi}). Un exemple d'espace 
$(\Omega, \mathcal{F}, \mathbb{P}, (\mathcal{F}_t)_{t \geq 0})$ 
satisfaisant la propri\'et\'e (NP) peut \^etre construit de la mani\`ere suivante:
\begin{itemize}
\item $\Omega$ est l'espace des fonctions continues de $\mathbb{R}_+$ and $\mathbb{R}^d$, ou
l'espace des fonctions c\`adl\`ag de $\mathbb{R}_+$ dans $\mathbb{R}^d$, pour un entier $d \geq 1$;
\item on d\'efinit $(\mathcal{F}^0_t)_{t \geq 0}$ comme la filtration canonique associ\'ee \`a 
$\Omega$, $\mathcal{F}^0$ comme la tribu engendr\'ee par $(\mathcal{F}^0_t)_{t \geq 0}$, et on 
fixe une mesure de probabilit\'e $\mathbb{P}^0$ sur $(\Omega, \mathcal{F}^0)$; 
\item on d\'efinit $(\Omega, \mathcal{F}, \mathbb{P}, (\mathcal{F}_t)_{t \geq 0})$ comme la 
plus petite extension possible de $(\Omega, \mathcal{F}^0, \mathbb{P}^0, (\mathcal{F}^0_t)_{t \geq 0})$ 
telle que pour tout $t \geq 0$, $\mathcal{F}_0$ contient toutes les parties de $\Omega$ incluses 
dans un ensemble $A \in \mathcal{F}^0_t$ tel que $\mathbb{P}^0 [A] = 0$. 
\end{itemize}
\noindent
Dans cette note, nous simplifions la preuve du Theor\`eme \ref{tout} donn\'ee dans \cite{NN1}. 
Cependent, la preuve de \cite{NN1} conserve son int\'er\^et dans la mesure o\`u elle est utilis\'ee 
pour d\'emontrer certaines propri\'etes importantes de la mesure $\mathcal{Q}$ (donn\'es dans \cite{NN1} et 
\cite{NN3}). Par ailleurs, nous prouvons une version du Th\'eor\`eme \ref{tout} valable pour des 
sous-martingales \`a temps discret. L'\'enonc\'e est le suivant:
\begin{theorem} \label{toutdiscret}
Soit $(X_n)_{n \geq 0}$ une sous-martingale \`a temps discret (en particulier $X_n$ est int\'egrable
pour tout $n \geq 0$), d\'efinie sur un espace de probabilit\'e filtr\'e 
$(\Omega, \mathcal{F}, \mathbb{P}, (\mathcal{F}_n)_{n \geq 0})$ qui satisfait une certaine propri\'et\'e 
technique pr\'cis\'ee plus bas.
On suppose que pour tout $n \geq 0$, $X_n = N_n + A_n$, o\`u 
$(N_n)_{n \geq 0}$ et $(A_n)_{n \geq 0}$ sont des processus 
satisfaisant les conditions suivantes:
\begin{itemize}
\item $(N_n)_{n \geq 0}$ est une martingale d\'efinie sur l'espace 
 $(\Omega, \mathcal{F}, \mathbb{P}, (\mathcal{F}_n)_{n \geq 0})$ 
\item $(A_n)_{n \geq 0}$ est un processus croissant, pr\'evisible (i.e. $A_n$ est $\mathcal{F}_{n-1}$-mesurable
pour tout $n \geq 1$), tel que $A_0 = 0$;
\item Pour tout $n \geq 0$, $(A_{n+1} - A_n)X_n = 0$, i.e. $A$ ne peut cro\^itre qu'apr\`es les z\'eros de
$X$. 
\end{itemize}
\noindent
Dans ces conditions, il existe une unique mesure $\sigma$-finie $\mathcal{Q}$, d\'efinie sur 
$(\Omega, \mathcal{F}, \mathbb{P})$ et telle que si $g = \sup\{n \geq 0, X_n = 0 \}$:
\begin{itemize}
\item $\mathcal{Q} [g = \infty] = 0$;
\item Pour tout $n \geq 1$, et pour toute variable al\'eatoire $F_n$ born\'ee, $\mathcal{F}_n$-mesurable,
$$\mathcal{Q} \left[ F_n \, \mathds{1}_{g < n} \right] = \mathbb{E}_{\mathbb{P}} \left[F_n X_n \right].$$
\end{itemize}
\noindent
\end{theorem}
\noindent
La preuve 
du Th\'eor\`me \ref{toutdiscret} est exactement similaire \`a la nouvelle preuve du Th\'eor\`eme \ref{tout}. 
\section{The continuous-time case}
The submartingales of class $(\Sigma)$ were introduced by Yor in \cite{Y}, and some of their main
properties are studied by Nikeghbali in \cite{N}. These submartingales are, in particular, related 
with some problems in mathematical finance (see \cite{MRY}, \cite{BeYor}, \cite{CNP}) and 
with Brownian penalisations (see \cite{NRY}). The class $(\Sigma)$ is defined as follows: 
\begin{definition}
Let $(\Omega,\mathcal{F},(\mathcal{F}_t)_{t \geq 0},\mathbb{P})$ be a filtered probability space. A 
nonnegative submartingale $(X_t)_{t \geq 0}$ is of class $(\Sigma)$, if it can be decomposed as
$X_t = N_t + A_t$ where $(N_t)_{t \geq 0}$ and $(A_t)_{t \geq 0}$ are $(\mathcal{F}_t)_{t \geq 0}$-adapted 
processes satisfying the following assumptions:
\begin{itemize}
\item $(N_t)_{t \geq 0}$ is a c\`adl\`ag  martingale;
\item $(A_t)_{t \geq 0}$ is a continuous increasing process, with $A_0 = 0$;
\item The measure $(dA_t)$ is carried by the set $\{t \geq 0, X_t = 0 \}$.
\end{itemize}
\end{definition}
\noindent
In \cite{NN1}, we prove that one can associate a $\sigma$-finite measure to any submartingale 
of class $(\Sigma)$, defined on a filtered probability space 
$(\Omega,\mathcal{F},(\mathcal{F}_t)_{t \geq 0},\mathbb{P})$ satisfying a technical condition, 
called property (NP). The statement is the following:
\begin{theorem} \label{all}
Let $(X_t)_{t \geq 0}$ be a submartingale of the class $(\Sigma)$ (in particular $X_t$ is integrable 
for all $t \geq 0$), defined on a filtered probability space 
$(\Omega, \mathcal{F}, \mathbb{P}, (\mathcal{F}_t)_{t \geq 0})$ which satisfies the property (NP).
In particular, $(\mathcal{F}_t)_{t \geq 0}$
 is right-continuous and $\mathcal{F}$ is 
the $\sigma$-algebra generated by $\mathcal{F}_t$, $t \geq 0$. 
Then, there exists a unique $\sigma$-finite measure $\mathcal{Q}$, defined 
on $(\Omega, \mathcal{F}, \mathbb{P})$, 
such that for $g:= \sup\{t \geq 0, X_t = 0 \}$:
\begin{itemize}
\item $\mathcal{Q} [g = \infty] = 0$;
\item For all $t \geq 0$, and for all $\mathcal{F}_t$-measurable, bounded random variables $F_t$,
$$\mathcal{Q} \left[ F_t \, \mathds{1}_{g \leq t} \right] = \mathbb{E}_{\mathbb{P}} \left[F_t X_t \right].$$
\end{itemize}
\noindent
\end{theorem} 
The property (NP) is precisely defined in \cite{NN2}, and is involved in particular because of the problems
encountered in the extension of compatible families of probability measures (these problems are also 
discussed by Bichteler in \cite{Bi}).  An example of space  
$(\Omega, \mathcal{F}, \mathbb{P}, (\mathcal{F}_t)_{t \geq 0})$ 
satisfying the property (NP) can be constructed in the following way:
\begin{itemize}
\item $\Omega$ is the space of continuous functions from $\mathbb{R}_+$ to $\mathbb{R}^d$, or
the space of c\`adl\`ag functions from $\mathbb{R}_+$ to $\mathbb{R}^d$, for an integer $d \geq 1$;
\item one defines $(\mathcal{F}^0_t)_{t \geq 0}$ as the canonical filtration associated with 
$\Omega$, $\mathcal{F}^0$ the $\sigma$-algebra generated by $(\mathcal{F}^0_t)_{t \geq 0}$, and 
one fixes a probability measure $\mathbb{P}^0$ on $(\Omega, \mathcal{F}^0)$; 
\item one defines $(\Omega, \mathcal{F}, \mathbb{P}, (\mathcal{F}_t)_{t \geq 0})$ as the smallest 
possible extension of $(\Omega, \mathcal{F}^0, \mathbb{P}^0, (\mathcal{F}^0_t)_{t \geq 0})$ 
such that for all $t \geq 0$, $\mathcal{F}_0$ contains all the subsets of $\Omega$ included in a set
 $A \in \mathcal{F}^0_t$ such that $\mathbb{P}^0 [A] = 0$. 
\end{itemize}
\noindent
In this note, we simplify the proof of Theorem \ref{all} given in \cite{NN1}. However, the proof in \cite{NN1}
remains interesting, since it is used to establish  some important properties of the measure 
$\mathcal{Q}$ (given in \cite{NN1} and \cite{NN3}).
Moreover, we are able to prove an analog of Theorem \ref{all} for a certain class of 
discrete-time submartingales. The precise statement in given in Section \ref{discrete}.
The new proof of Theorem \ref{all} and the proof of the discrete-time result are given in Section \ref{proofs}:
these two proofs are very similar. 
\section{The discrete-time case} \label{discrete}
The discrete-time result is very similar to Theorem \ref{all}. However, one needs to define a discrete 
version of the property (NP). This can be done in the same way as for the continuous 
 case. We first state the following definition:
\begin{definition} \label{P}
Let $(\Omega, \mathcal{F}, (\mathcal{F}_n)_{n \geq 0})$ be a filtered measurable space, such that
$\mathcal{F}$ is the $\sigma$-algebra generated by $\mathcal{F}_n$, 
$n \geq 0$: $\mathcal{F}=\bigvee_{n\geq0}\mathcal{F}_n$. We shall say that the
 property (Pd) holds if and only if $(\mathcal{F}_n)_{n \geq 0}$ 
enjoys the following conditions: 
\begin{itemize}
\item For all $n \geq 0$, $\mathcal{F}_n$ is generated by a countable number of sets;
\item For all $n \geq 0$, there exists a Polish space $\Omega_n$, and a surjective map 
 $\pi_n$ from $\Omega$ to $\Omega_n$, such that $\mathcal{F}_n$ is the $\sigma$-algebra of the inverse
 images, by $\pi_n$, of Borel sets in $\Omega_n$, and such that for all $B \in \mathcal{F}_n$, 
 $\omega \in \Omega$, $\pi_n (\omega) \in \pi_n(B)$ implies $\omega \in B$;
\item If $(\omega_n)_{n \geq 0}$ is a sequence of elements of $\Omega$, such that for all $N \geq 0$,
$$\bigcap_{n = 0}^{N} A_n (\omega_n) \neq \emptyset,$$
where $A_n (\omega_n)$ is the intersection of the sets in $\mathcal{F}_n$ containing $\omega_n$, 
then:
$$\bigcap_{n = 0}^{\infty} A_n (\omega_n) \neq \emptyset.$$
\end{itemize}
\end{definition}
\noindent
A fundamental example of space satisfying the property (Pd) is obtained by taking
$\Omega = X^{\mathbb{N}}$ where $X$ is a Polish space, and 
$\mathcal{F}=\bigvee_{n\geq0}\mathcal{F}_n$, where for $n \geq 0$, $\mathcal{F}_n$ is
 the $\sigma$-algebra of the inverse images, by the projection on the $n$ first 
coordinates, of Borel sets in $X^n$. 
Another example can be constructed as follows: 
\begin{itemize}
\item one defines $\Omega$ as the space $\mathcal{C}(\mathbb{R}_+,\mathbb{R}^d)$ of
 continuous functions from 
	$\mathbb{R}_+$ to $\mathbb{R}^d$, or as the space $\mathcal{D}(\mathbb{R}_+,\mathbb{R}^d)$ 
of c\`adl\`ag functions
 from $\mathbb{R}_+$ 
	to $\mathbb{R}^d$ (for some $d \geq 1$);
\item for $n \geq 0$, one defines $(\mathcal{F}_n)_{n \geq 0}$
as the natural filtration of the canonical process,
and $\mathcal{F}=\bigvee_{n\geq0}\mathcal{F}_n$. 
\end{itemize}
The interest of property (Pd) lies in the following result:
\begin{proposition} \label{extension}
Let $(\Omega, \mathcal{F}, (\mathcal{F}_n)_{n \geq 0})$ be a filtered measurable space satisfying the 
property (Pd), and let, for $n \geq 0$, $(\mathbb{Q}_n)$ be a family of probability measures on $(\Omega, \mathcal{F}_n)$, 
such that for all $n \geq m \geq 0$, $\mathbb{Q}_m$ is the restriction of $\mathbb{Q}_n$ to $\mathcal{F}_m$.
Then, there exists a unique measure $\mathbb{Q}$ on $(\Omega, \mathcal{F})$ such that for all $n \geq 0$,
its restriction to $\mathcal{F}_n$ is equal to $\mathbb{Q}_n$. 
\end{proposition}
\noindent
Moreover, one can combine the property (Pd) with the natural augmentation, defined by the following 
result:
 \begin{proposition}
Let $(\Omega, \mathcal{F}, (\mathcal{F}_n)_{n \geq 0}, \mathbb{P})$ be a filtered probability space, 
and let $\mathcal{N}_0$ be the family of the subsets $A$ of $\Omega$ which satisfy the following property:
there exists $k \geq 0$ and $B \in \mathcal{F}_k$ such that $A \subset B$ and $\mathbb{P}[B]= 0$. 
One defines $\widetilde{\mathcal{F}}$ as the $\sigma$-algebra generated by $\mathcal{F}$ and 
$\mathcal{N}_0$, and for $n \geq 0$, $\widetilde{\mathcal{F}}_n$ as the $\sigma$-algebra generated 
by $\mathcal{F}_n$ and $\mathcal{N}_0$. Then, there exists a unique extension $\widetilde{\mathbb{P}}$
of $\mathbb{P}$, defined on $(\Omega, \widetilde{\mathcal{F}})$. The space 
$(\Omega, \widetilde{\mathcal{F}}, (\widetilde{\mathcal{F}}_n)_{n \geq 0}, \widetilde{\mathbb{P}})$ 
is called the natural augmentation of the space
$(\Omega, \mathcal{F}, (\mathcal{F}_n)_{n \geq 0}, \mathbb{P})$.
\end{proposition}
\noindent
Note that in discrete case, there does not exist a notion of right-continuous filtration. That 
is why the definition of the natural augmentation is slightly different from the definition in the
continuous time case. The extension of measures is still available when one takes the natural augmentation 
of a space satisfying the property (Pd). More precisely, one has the following result:
\begin{proposition} \label{extensionaugmentation}
Let $(\Omega, \mathcal{F}, (\mathcal{F}_n)_{n \geq 0}, \mathbb{P})$ be the 
 natural augmentation of a filtered probability space satisfying the property (Pd). 
 Then, for all coherent families of probability measures
$(\mathbb{Q}_n)_{n \geq 0}$, such that $\mathbb{Q}_n$ is defined 
on $\mathcal{F}_n$,
and is absolutely continuous with respect to the restriction
of $\mathbb{P}$ to $\mathcal{F}_n$, there exists a unique probability measure
$\mathbb{Q}$
on $\mathcal{F}$ which coincides with $\mathbb{Q}_n$ on $\mathcal{F}_n$ 
for all $n \geq 0$. 
\end{proposition}
\noindent
The proofs of all these results are essentially given in \cite{NN2} (the only change is the replacement
of continuous filtrations by discrete filtrations). 
 The method used in 
our proof of Proposition \ref{extension} comes from Stroock and Varadhan (see \cite{SV}). 
Note that the condition (Pd) is 
not new and is essentially given by Parthasarathy in \cite{Parth}, p. 141. 
One can now state the discrete analog of Theorem \ref{all}
\begin{theorem} \label{alldiscrete}
Let $(\Omega, \mathcal{F}, \mathbb{P}, (\mathcal{F}_n)_{n \geq 0})$ be a filtered probability space 
satisfying the property (Pd), or its natural augmentation. 
One considers a discrete-time submartingale $(X_n)_{n \geq 0}$ (in particular, $X_n$ is integrable
for all $n \geq 0$), defined on this space, and one supposes that for
 all $n \geq 0$, $X_n = N_n + A_n$, where the processes 
$(N_n)_{n \geq 0}$ and $(A_n)_{n \geq 0}$ satisfy the following conditions:
\begin{itemize}
\item $(N_n)_{n \geq 0}$ is a martingale defined on the space
 $(\Omega, \mathcal{F}, \mathbb{P}, (\mathcal{F}_n)_{n \geq 0})$ 
\item $(A_n)_{n \geq 0}$ is an increasing, predictable process ($A_n$ is $\mathcal{F}_{n-1}$-measurable
for all $n \geq 1$), such that $A_0 = 0$;
\item For all $n \geq 0$, $(A_{n+1} - A_n)X_n = 0$, i.e. $A$ can only increase after the zeros of 
$X$. 
\end{itemize}
\noindent
Under these assumptions, there exists a unique $\sigma$-finite measure $\mathcal{Q}$, 
defined on 
$(\Omega, \mathcal{F}, \mathbb{P})$ and such that for $g := \sup\{n \geq 0, X_n = 0 \}$:
\begin{itemize}
\item $\mathcal{Q} [g = \infty] = 0$;
\item For all $n \geq 1$, and for all bounded, $\mathcal{F}_n$-measurable random variable $F_n$, 
$$\mathcal{Q} \left[ F_n \, \mathds{1}_{g < n} \right] = \mathbb{E}_{\mathbb{P}} \left[F_n X_n \right].$$
\end{itemize}
\noindent
\end{theorem}
\noindent
The new proof of Theorem \ref{all} and the proof and Theorem \ref{alldiscrete} are given in Section 
\ref{proofs}.
\section{Proof of Theorems \ref{all} and \ref{alldiscrete}} \label{proofs}
Let us first proof Theorem \ref{all}. For $t \geq 0$, we define $d_t$ as the smallest hitting time of 
zero strictly after time $t$:
$$d_t := \inf \{s > t, X_s = 0\}.$$
By the d\'ebut theorem, available for the spaces satisfying the property (NP) (as proved in \cite{NN2}), 
$d_t$ is a stopping time. One can prove that $(X_{s \wedge d_t})_{s \geq t}$ is a martingale
with repect to the filtration $(\mathcal{F}_s)_{s \geq t}$ and the probability measure $\mathbb{P}$. 
Indeed, for all $s \geq t$, $$X_{s \wedge d_t} = N_{s \wedge d_t} + A_{s \wedge d_t} = N_{s \wedge d_t} + A_t,$$
where $(N_{s \wedge d_t})_{s \geq 0}$ is a martingale by the stopping theorem, and $A_t$
 is $\mathcal{F}_s$-measurable and 
independent of $s$. Since $(\Omega, \mathcal{F}, \mathbb{P}, (\mathcal{F}_t)_{t \geq 0})$ satisfies 
the property (NP), there exists a finite measure $\mathcal{Q}^{(t)}$ on $(\Omega, \mathcal{F})$ such that
for all $s \geq t$, and for all events $\Lambda_s \in \mathcal{F}_s$:
$$ \mathcal{Q}^{(t)} [\Lambda_s]  = \mathbb{E}_{\mathbb{P}} [\mathds{1}_{\Lambda_s} X_{s \wedge d_t} ].$$
Once $\mathcal{Q}^{(t)}$ is constructed, let us observe that under this measure, $g \leq t$ almost everywhere. 
Indeed, for all $s \geq t$:
$$\mathcal{Q}^{(t)} [d_t \leq s] = \mathbb{E}_{\mathbb{P}} [\mathds{1}_{d_t \leq s} X_{d_t}]= 0,$$
since for $d_t < \infty$, $X_{d_t}$ vanishes by right-continuity of $X$. 
Hence, $\mathcal{Q}^{(t)}$-almost everywhere, $d_t = \infty$, or equivalently, $g \leq t$. 
Let us now prove that for $u \geq t \geq 0$, $\mathcal{Q}^{(t)}$ is the restriction 
of $\mathcal{Q}^{(u)}$ to the set $\{g \leq t\}$. Since $\mathcal{Q}^{(u)}$
almost-everywhere, $X$ does not vanish strictly after time $u$, one has for all $s \geq u$, and for
 all events $\Lambda_s \in \mathcal{F}_s$:
\begin{align*}
\mathcal{Q}^{(u)}[g \leq t, \Lambda_s] & = 
\mathcal{Q}^{(u)} [d_t = \infty, \Lambda_s] \\ & = 
\mathcal{Q}^{(u)} [d_t > s, \Lambda_s] \\ & = 
\mathbb{E}_{\mathbb{P}} [\mathds{1}_{d_t > s, \Lambda_s} X_{d_u \wedge s}] \\
& =\mathbb{E}_{\mathbb{P}} [\mathds{1}_{d_t > s, \Lambda_s} X_{s}] \\ & = 
\mathbb{E}_{\mathbb{P}} [\mathds{1}_{\Lambda_s} X_{d_t \wedge s}] \\
& = \mathcal{Q}^{(t)} [\Lambda_s]
\end{align*}
\noindent
Since $\mathcal{Q}^{(t)}$ is the restriction of $\mathcal{Q}^{(u)}$ to the set $\{g \leq t\}$, 
the family of measures $(\mathcal{Q}^{(t)})_{t \geq 0}$ is increasing. One can then define, for 
all events $\Lambda \in \mathcal{F}$:
$$\mathcal{Q} [\Lambda] = \underset{t \rightarrow \infty}{\lim} \mathcal{Q}^{(t)} [\Lambda],$$
where the limit is increasing. It is easy to check that $\mathcal{Q}$ is countably additive, hence, it is
a $\sigma$-finite measure. Moreover, $g$ is finite $\mathcal{Q}^{(t)}$-almost everywhere for all $t \geq 0$, 
and hence, $\mathcal{Q}$-almost everywhere. 
Now, for all $u \geq t \geq 0$, $\Lambda_t \in \mathcal{F}_t$:
$$\mathcal{Q}^{(u)} [g \leq t, \Lambda_t] = \mathcal{Q}^{(t)} [\Lambda_t] = \mathbb{E}_{\mathbb{P}}
[\mathds{1}_{\Lambda_t} X_t],$$
and by taking $u \rightarrow \infty$:
$$\mathcal{Q} [g \leq t, \Lambda_t] =\mathbb{E}_{\mathbb{P}} [\mathds{1}_{\Lambda_t} X_t].$$
We have now proved the existence part of Theorem \ref{all}. The uniqueness part is very easily 
shown in \cite{NN1}. Let us now sketch the proof of Theorem \ref{alldiscrete}, which is similar 
to the proof of Theorem \ref{all} given just above. 
For $n \geq 0$, 
let us define $d_n$ as the first hitting time of zero strictly after time $n$:
the process $(X_{p \wedge d_n})_{p \geq n+1}$ is a martingale with respect to 
the filtration $(\mathcal{F}_p)_{p \geq n+1}$. By Propositions \ref{extension} and 
\ref{extensionaugmentation}, one deduces that there exists 
a finite measure $\mathcal{Q}^{(n)}$ such that its density with respect to $\mathbb{P}$, after 
restriction to $\mathcal{F}_p$, is equal to $X_{p \wedge d_n}$, for all $p \geq n+1$. 
Moreover, one has $g \leq n$, $\mathcal{Q}^{(n)}$-almost everywhere, and for $n \geq m$, 
$\mathcal{Q}^{(m)}$ is the restriction of $\mathcal{Q}^{(n)}$ to the set $g \leq m$. 
Hence, one can define $\mathcal{Q}$ as the increasing limit of $\mathcal{Q}^{(n)}$ when
$n$ goes to infinity. One obtains, for all $n \geq 1$, and for all events $\Lambda_n \in \mathcal{F}_n$:
$$\mathcal{Q} [g < n, \Lambda_n] = \mathcal{Q}^{(n-1)} [\Lambda_n] = 
\mathbb{E}_{\mathbb{P}} [\mathds{1}_{\Lambda_n} X_n].$$
The uniqueness is proved as in the continuous-time case (see \cite{NN1}).

\providecommand{\bysame}{\leavevmode\hbox to3em{\hrulefill}\thinspace}
\providecommand{\MR}{\relax\ifhmode\unskip\space\fi MR }
\providecommand{\MRhref}[2]{%
  \href{http://www.ams.org/mathscinet-getitem?mr=#1}{#2}
}
\providecommand{\href}[2]{#2}

\end{document}